\documentclass[11pt]{article} 
\usepackage{latexsym} 
\usepackage{amsmath}
\usepackage{amssymb} 
\usepackage{amsfonts}
\usepackage{amscd} 
\usepackage{graphicx}
\usepackage{layout}
\begin{document} 
\newtheorem{Th}{Theorem}[section]
\newtheorem{Cor}{Corollary}[section]
\newtheorem{Prop}{Proposition}[section]
\newtheorem{Lem}{Lemma}[section]
\newtheorem{Def}{Definition}[section]
\newtheorem{Rem}{Remark}[section]
\newtheorem{Ex}{Example}[section]
\newtheorem{stw}{Proposition}[section]


\newcommand{\bet}{\begin{Th}}
\newcommand{\ent}{\stepcounter{Cor}
   \stepcounter{Prop}\stepcounter{Lem}\stepcounter{Def}
   \stepcounter{Rem}\stepcounter{Ex}\end{Th}}


\newcommand{\bec}{\begin{Cor}}
\newcommand{\enc}{\stepcounter{Th}
   \stepcounter{Prop}\stepcounter{Lem}\stepcounter{Def}
   \stepcounter{Rem}\stepcounter{Ex}\end{Cor}}
\newcommand{\bep}{\begin{Prop}}
\newcommand{\enp}{\stepcounter{Th}
   \stepcounter{Cor}\stepcounter{Lem}\stepcounter{Def}
   \stepcounter{Rem}\stepcounter{Ex}\end{Prop}}
\newcommand{\bel}{\begin{Lem}}
\newcommand{\enl}{\stepcounter{Th}
   \stepcounter{Cor}\stepcounter{Prop}\stepcounter{Def}
   \stepcounter{Rem}\stepcounter{Ex}\end{Lem}}
\newcommand{\bef}{\begin{Def}}
\newcommand{\enf}{\stepcounter{Th}
   \stepcounter{Cor}\stepcounter{Prop}\stepcounter{Lem}
   \stepcounter{Rem}\stepcounter{Ex}\end{Def}}
\newcommand{\ber}{\begin{Rem}}
\newcommand{\enr}{
   \stepcounter{Th}\stepcounter{Cor}\stepcounter{Prop}
   \stepcounter{Lem}\stepcounter{Def}\stepcounter{Ex}\end{Rem}}
\newcommand{\bee}{\begin{Ex}}
\newcommand{\ene}{
   \stepcounter{Th}\stepcounter{Cor}\stepcounter{Prop}
   \stepcounter{Lem}\stepcounter{Def}\stepcounter{Rem}\end{Ex}}
\newcommand{\Proof}{\noindent{\it Proof\,}:\ }
\newcommand{\beP}{\Proof}
\newcommand{\enP}{\hfill $\Box$ \par\vspace{5truemm}}

\newcommand{\EE}{\mathbf{E}}
\newcommand{\QQ}{\mathbf{Q}}
\newcommand{\R}{\mathbf{R}}
\newcommand{\C}{\mathbf{C}}
\newcommand{\ZZ}{\mathbf{Z}}
\newcommand{\KK}{\mathbf{K}}
\newcommand{\NN}{\mathbf{N}}
\newcommand{\PP}{\mathbf{P}}
\newcommand{\HH}{\mathbf{H}}
\newcommand{\uuu}{\boldsymbol{u}}
\newcommand{\xxx}{\boldsymbol{x}}
\newcommand{\aaa}{\boldsymbol{a}}
\newcommand{\bbb}{\boldsymbol{b}}
\newcommand{\AAA}{\mathbf{A}}
\newcommand{\BBB}{\mathbf{B}}
\newcommand{\ccc}{\boldsymbol{c}}
\newcommand{\iii}{\boldsymbol{i}}
\newcommand{\jjj}{\boldsymbol{j}}
\newcommand{\kkk}{\boldsymbol{k}}
\newcommand{\rrr}{\boldsymbol{r}}
\newcommand{\FFF}{\boldsymbol{F}}
\newcommand{\yyy}{\boldsymbol{y}}
\newcommand{\ppp}{\boldsymbol{p}}
\newcommand{\qqq}{\boldsymbol{q}}
\newcommand{\nnn}{\boldsymbol{n}}
\newcommand{\vvv}{\boldsymbol{v}}
\newcommand{\eee}{\boldsymbol{e}}
\newcommand{\fff}{\boldsymbol{f}}
\newcommand{\www}{\boldsymbol{w}}
\newcommand{\0}{\boldsymbol{0}}
\newcommand{\lon}{\longrightarrow}
\newcommand{\ga}{\gamma}
\newcommand{\pa}{\partial}
\newcommand{\QED}{\hfill $\Box$}
\newcommand{\id}{{\mbox {\rm id}}}
\newcommand{\Ker}{{\mbox {\rm Ker}}}
\newcommand{\grad}{{\mbox {\rm grad}}}
\newcommand{\ind}{{\mbox {\rm ind}}}
\newcommand{\rot}{{\mbox {\rm rot}}}
\newcommand{\diver}{{\mbox {\rm div}}}
\newcommand{\Gr}{{\mbox {\rm Gr}}}
\newcommand{\LG}{{\mbox {\rm LG}}}
\newcommand{\Diff}{{\mbox {\rm Diff}}}
\newcommand{\Symp}{{\mbox {\rm Symp}}}
\newcommand{\Ct}{{\mbox {\rm Ct}}}
\newcommand{\Uns}{{\mbox {\rm Uns}}}
\newcommand{\rank}{{\mbox {\rm rank}}}
\newcommand{\sign}{{\mbox {\rm sign}}}
\newcommand{\Spin}{{\mbox {\rm Spin}}}
\newcommand{\Sp}{{\mbox {\rm sp}}}
\newcommand{\Int}{{\mbox {\rm Int}}}
\newcommand{\Hom}{{\mbox {\rm Hom}}}
\newcommand{\Tan}{{\mbox {\rm Tan}}}
\newcommand{\codim}{{\mbox {\rm codim}}}
\newcommand{\ord}{{\mbox {\rm ord}}}
\newcommand{\Iso}{{\mbox {\rm Iso}}}
\newcommand{\corank}{{\mbox {\rm corank}}}
\def\mod{{\mbox {\rm mod}}}
\newcommand{\pt}{{\mbox {\rm pt}}}
\newcommand{\qed}{\hfill $\Box$ \par}
\newcommand{\spe}{\vspace{0.4truecm}}
\newcommand{\ad}{{\mbox{\rm ad}}}

\newcommand{\dint}[2]{{\displaystyle\int}_{{\hspace{-1.9truemm}}{#1}}^{#2}}

%
\newenvironment{FRAME}{\begin{trivlist}\item[]
	\hrule
	\hbox to \linewidth\bgroup
		\advance\linewidth by -10pt
		\hsize=\linewidth
		\vrule\hfill
		\vbox\bgroup
			\vskip5pt
			\def\thempfootnote{\arabic{mpfootnote}}
			\begin{minipage}{\linewidth}}{%
			\end{minipage}\vskip5pt
		\egroup\hfill\vrule
	\egroup\hrule
	\end{trivlist}}

\title{
Singularities of tangent surfaces to directed curves
} 

\author{G. Ishikawa\thanks{This work was supported by JSPS KAKENHI No.15H03615 and No.15K13431.} \ 
and T. Yamashita}

\renewcommand{\thefootnote}{\fnsymbol{footnote}}

\date{ }

\maketitle

\begin{abstract} 
A directed curve is a possibly singular 
curve with well-defined tangent lines along the curve. 
Then the tangent surface to a directed curve is naturally defined 
as the ruled surface by tangent geodesics to the curve, 
whenever any affine connection is endowed with the ambient space. 
In this paper the local diffeomorphism classification is 
completed for generic directed curves. Then it turns out that 
the swallowtails and open swallowtails appear generically 
for the classification on singularities of tangent surfaces. 
\end{abstract} 

\section{Introduction}
\label{Introduction}

Given a space curve, 
the ruled surface by its tangent lines is called a {\it tangent surface} or a {\it tangent developable} to the curve. 
Tangent surfaces appear in various geometric problems and applications (see for instance \cite{BG}\cite{IMT0}). 
Even if the space curve is regular, its tangent surface has singularities 
at least along the original curve, so called \lq\lq the curve of regression\rq\rq. 

Let $M$ be a general (semi-)Riemannian manifold, or more generally, 
a manifold $M$ with an affine connection $\nabla$, of dimension $m \geq 3$, 
and let $\gamma : I \to M$ any regular curve in $M$. 
If we replace tangent lines by \lq\lq tangent geodesics\rq\rq \, in the definition of tangent surface, 
then we have the definition of the $\nabla$-tangent surface 
$\nabla$-$\Tan(\gamma) : (I \times \R, I \times \{ 0\}) \to M$ as a map-germ along $I\times\{ 0\}$. 

Ordinarily we try to classify certain generic singularities in a {\it specific} space, say, 
in the Euclidian spaces, in the space forms, and so on. 
If we treat arbitrary spaces, it would become hopeless to classify singularities of tangent surfaces that appear far away. 
However, it is possible to find a local classification theorem which holds in general spaces. 
In the previous paper \cite{IY}, actually we have shown the following result on 
the singularities of $\nabla$-tangent surfaces to generic curves for arbitrary affine connection $\nabla$: 

\bet
\label{genericity-theorem-1} 
{\rm (\cite{IY})}
The singularities of the $\nabla$-tangent surface 
to a generic immersed curve in $M$ on a neighbourhood of the curve 
are only the cuspidal edges and the folded umbrellas if $m = 3$, 
and the embedded cuspidal edges if $m \geq 4$. 
\ent

The above theorem provides a rare but an ultimate {\it local} classification of singularities 
associated with {\it generic} immersed curves in {\it general} spaces. 
The explanation on singularities is coming later soon. 

Now regarding the definition of general tangent surfaces, it seems to be very natural to consider the 
genericity in the space of curves, not only for regular (immersed) curves, 
but also for all singular curves with well-defined tangent directions, called {\it directed curves}, and 
to classify singularities of tangent surfaces for curves which is generic in such a class. 
In fact, as we show in this paper, it is possible and we have the following general result: 

\bet
\label{genericity-theorem} {\rm (Singularities of tangent surfaces to generic directed curves.) }
Let $\nabla$ be any affine connection on a manifold $M$ of dimension $m \geq 3$. 
The singularities of the $\nabla$-tangent surface to a generic directed 
curve in $M$ on a neighbourhood of the curve 
are only the cuspidal edges, the folded umbrellas and the swallowtails if $m = 3$,  
and the embedded cuspidal edges and open swallowtails if $m \geq 4$. 
\ent

The genericity is exactly given using 
Whitney $C^\infty$ topology on an appropriate space of curves (see Proposition \ref{genericity2}).

\

A map-germ $f : (\R^2, p) \to M$ is locally {\it diffeomorphic} at $p$ to another map-germ 
$g : (\R^2, p') \to M'$
if there exist diffeomorphism-germs 
$\sigma : (\R^2, p) \to (\R^2, p')$ and $\tau : (M, f(p)) \to (M', g(p'))$ such that 
$\tau\circ f = g\circ \sigma : (\R^2, p) \to (M', g(p'))$. 

The {\it cuspidal edge} is defined by the map-germ $(\R^2, 0) \to (\R^m, 0)$, $m \geq 3$, 
$$
(t, s) \mapsto (t + s, \ t^2 + 2st, \ t^3 + 3st^2, \ 0, \ \dots, \ 0), 
$$
which is diffeomorphic to $(u, w) \mapsto (u, w^2, w^3, 0, \dots, 0)$. 
The cuspidal edge singularities are originally defined only in the three dimensional space. 
Here we are generalizing the notion of the cuspidal edge in higher dimensional space. In Theorem \ref{genericity-theorem}, we emphasize it by writing \lq\lq embedded" cuspidal edge. In what follows, we call it just cuspidal edge for simplicity even in the case $m \geq 4$. 
The {\it folded umbrella} (or the {\it cuspidal cross cap}) is defined by the map-germ 
$(\R^2, 0) \to (\R^3, 0)$, 
$$
(t, s) \mapsto (t + s, \ t^2 + 2st, \ t^4 + 4st^3), 
$$
which is diffeomorphic to $(u, t) \mapsto (u, t^2 + ut, t^4 + {\textstyle \frac{2}{3}}ut^3)$. 
The {\it swallowtail} is defined by the map-germ $(\R^2, 0) \to (\R^3, 0)$ 
$$
(t, s) \mapsto (t^2 + s, \ t^3 + {\textstyle \frac{3}{2}}st, \ t^4 + 2st^2), 
$$
which is diffeomorphic to $(u, t) \mapsto (u, t^3 + ut, t^4 + \frac{2}{3}ut^2)$. 
The {\it open swallowtail} is defined by the map-germ $(\R^2, 0) \to (\R^m, 0)$, $m \geq 4$, 
$$
(t, s) \mapsto (t^2 + s, \ t^3 + {\textstyle \frac{3}{2}}st, \ t^4 + 2st^2, \ t^5 + {\textstyle \frac{5}{2}}st^3, \ 0, \ \dots, \ 0), 
$$
which is diffeomorphic to $(u, t) \mapsto (u, t^3 + ut, t^4 + \frac{2}{3}ut^2, t^5 + \frac{5}{9}ut^3, 0, \dots, 0)$.
The open swallowtail singularity was introduced by Arnol'd (see \cite{Arnold}) as a singularity of Lagrangian varieties in symplectic geometry. Here we abstract its diffeomorphism class as the singularity of tangent surfaces (see \cite{Givental}\cite{Ishikawa4}). 



Swallowtails and open swallowtails appear as singularities of tangent surfaces to {\it singular} curves. 
It is observed that (open) swallowtails are destroyed by some perturbations of the original curves which 
induce big changes of their tangent directions, 
and however that they survive by any small perturbations which induce small changes of their tangent directions
of the singular but directed curves. 

\

Let $\gamma : I \to M$ be any curve which is not necessarily a geodesic nor 
an immersed curve. 
The first derivative $(\nabla\gamma)(t)$ means just the velocity 
vector field $\gamma'(t)$. The second derivative $(\nabla^2\gamma)(t)$ is defined, in terms of 
covariant derivative along the curve $\gamma$, by 
$$
(\nabla^2\gamma)(t) := \nabla_{\pa/\pa t}^\gamma (\nabla\gamma)(t). 
$$
Note that $\gamma$ is a $\nabla$-geodesic if and only if $\nabla^2 \gamma = 0$. 
In general, we define $k$-th covariant derivative of $\gamma$ inductively by 
$$
(\nabla^k\gamma)(t) := \nabla_{\pa/\pa t}^\gamma (\nabla^{k-1}\gamma)(t), \ (k \geq 2). 
$$
Then we have: 

\bet
\label{characterization-theorem}
{\rm (Characterization.)} 
Let $\nabla$ be a torsion free affine connection on a manifold $M$. 
Let $\gamma : I \to M$ be a $C^\infty$ curve from an open interval $I$. 

{\rm (1)}  Let $\dim(M) = 3$. 
If $(\nabla\gamma)(t_0), (\nabla^2\gamma)(t_0), (\nabla^3\gamma)(t_0)$ are linearly independent, 
then the $\nabla$-tangent surface 
$\nabla{\mbox{\rm -}}\Tan(\gamma)$ is locally diffeomorphic to the cuspidal edge at $(t_0, 0) \in I \times \R$. 
If $(\nabla\gamma)(t_0), (\nabla^2\gamma)(t_0), (\nabla^3\gamma)(t_0)$ are linearly dependent, 
and $(\nabla\gamma)(t_0), (\nabla^2\gamma)(t_0), (\nabla^4\gamma)(t_0)$ are linearly independent, then
$\nabla{\mbox{\rm -}}\Tan(\gamma)$ is 
locally diffeomorphic to the folded umbrella at $(t_0, 0) \in I \times \R$. 
If $(\nabla\gamma)(t_0) = 0$ and $(\nabla^2\gamma)(t_0), (\nabla^3\gamma)(t_0), (\nabla^4\gamma)(t_0)$ 
are linearly independent, then 
$\nabla{\mbox{\rm -}}\Tan(\gamma)$ is 
locally diffeomorphic to the swallowtail at $(t_0, 0) \in I \times \R$. 

{\rm (2)} Let $\dim(M) \geq 4$. 
If $(\nabla\gamma)(t_0), (\nabla^2\gamma)(t_0), (\nabla^3\gamma)(t_0)$ are linearly independent, 
then the $\nabla$-tangent surface $\nabla{\mbox{\rm -}}\Tan(\gamma)$ is locally diffeomorphic to the cuspidal edge at $(t_0, 0) \in I \times \R$. 
If $(\nabla\gamma)(t_0) = 0$ and 
$(\nabla^2\gamma)(t_0), (\nabla^3\gamma)(t_0), (\nabla^4\gamma)(t_0), (\nabla^5\gamma)(t_0)$ 
are linearly independent, then 
$\nabla{\mbox{\rm -}}\Tan(\gamma)$ is 
locally diffeomorphic to the open swallowtail at $(t_0, 0) \in I \times \R$. 
\ent

Some of characterizations in Theorem \ref{characterization-theorem} have been shown already in \cite{IY}.  

The intrinsic characterizations of singularities found in \cite{KRSUY}\cite{FSUY} 
are useful for our treatment of singularities in general ambient spaces. 
We apply to non-flat projective geometry the characterizations and their some generalization via the notion of openings introduced by the first author (\cite{Ishikawa4}, see also \cite{Ishikawa2}). 

\

In \S \ref{Tangent surface to directed curve} we introduce the notion of directed curves and 
define their tangent surfaces. 
We recall the criteria of singularities in 
\S \ref{Swallowtail and open swallowtail} and prove Theorem \ref{characterization-theorem}. 
In \S \ref{Proof of the genericity theorem} we study on perturbations of directed curves 
and prove Theorem \ref{genericity-theorem}. 

\

In this paper all manifolds and mappings are assumed to be of class $C^\infty$ unless otherwise stated. 

This paper is a second half of the unpublished paper \cite{IY-pre} which is divided into two shorter papers,  the paper \cite{IY} and the present paper. We utilize in the present paper, as the sequel of \cite{IY}, several detailed calculations performed in \cite{IY}.

\section{Directed curves and their tangent surfaces}
\label{Tangent surface to directed curve}

Let $PTM = \Gr(1, TM)$ denote the projective tangent bundle over the manifold $M$, 
and $\pi : PTM \to M$ the natural projection. The fibre of $\pi$ over $x \in M$ is the 
projective space $P(T_xM)$ of dimension $m-1$. 

A curve $\gamma : I \to M$ from an open interval $I$, which is not necessarily an immersion, is called {\it directed} 
if there assigned a $C^\infty$ lifting $\widetilde{\gamma} : I \to PTM$ of $\gamma$ for $\pi$ which 
satisfies the integrality condition $\gamma'(t) \in \widetilde{\gamma}(t) \subset T_{\gamma(t)}M$ for any $t \in I$. 
Here $\widetilde{\gamma}(t) \in P(T_{\gamma(t)}M)$ is regarded as a one-dimensional 
linear subspace of $T_{\gamma(t)}M$. 
Then we regard the direction $\widetilde{\gamma}(t_0)$ is assigned to 
each point $\gamma(t_0)$ on $\gamma$. Note that 
if $\gamma'(t_0) \not= 0$, then $\widetilde{\gamma}(t_0)$ is uniquely determined by 
the tangent line $\langle \gamma'(t_0)\rangle_{\R} \subset T_{\gamma(t_0)}M$. 
The notion of directed curves is nothing but the 
notion of frontal maps introduced in \cite{IY} in the case $n = 1$ 
with assignment of an integral lifting when the immersion locus of $\gamma$ is dense in $I$. 

Let $\gamma :  I \to M$ be a directed curve and $\widetilde{\gamma}$ its integral lifting. 
Then there exists a $C^\infty$ frame $u : I \to TM$ of $\widetilde{\gamma}$ which satisfies
$\widetilde{\gamma}(t) = \langle u(t)\rangle_{\R}, u(t) \not= 0$ for any $t \in I$. 
Note that there exists a unique function $a(t)$ such that 
$\gamma'(t) = a(t)u(t)$. 
Then define the $\nabla$-tangent surface $f = \nabla{\mbox{\rm -}}\Tan(\gamma) : V (\subset I \times \R) \to M$ by 
$$
f(t, s) := \varphi(\gamma(t), u(t), s), 
$$
using the family of $\nabla$-geodesics $\varphi = \varphi(x, v, s)$ and a frame $u(t)$. 
Here $\varphi(x, v, s)$ gives the $\nabla$-geodesic parametrized by the parameter $s$ 
through $x$ with the velocity vector $v$ at $s = 0$, $\varphi(x, v, 0) = x$ and $\frac{\pa \varphi}{\pa s}(x, v, 0) = v$. 
In \cite{IY}, the $\nabla$-tangent surface for an immersed curve $\gamma$ 
was defined by the frame $u(t) = \gamma'(t)$ and studied with the detail analysis of 
$\nabla$-geodesics $\varphi = \varphi(x, v, s)$. 

\bel
If the immersion locus of a directed curve $\gamma : I \to M$ is dense in $I$, 
then the integral lifting $\widetilde{\gamma}$ is uniquely determined. 
The right equivalence class of the germ of $\nabla{\mbox{\rm -}}\Tan(\gamma) : (I \times \R, I \times\{ 0\}) \to M$ 
for a directed curve $\gamma$ 
is independent of the choice of the frame $u$. 
\enl

\Proof
The first half is clear because $\widetilde{\gamma}$ is $C^\infty$, so is continuous. 
The second half is achieved by the diffeomorphism $(t, s) \to (t, c(t)s)$ for 
another choice $c(t)u(t), c(t) \not= 0$. 
\QED

\

In \cite{IY} we have introduced the notions of frontals and non-degenerate singular points of frontals 
(\S 3 of \cite{IY}). Using those notions we have the following result: 

\bel
\label{non-degenerate-frontal2}
Let $\gamma : I \to M$ be a $C^\infty$ curve, $t_0 \in I$, and $k \geq 1$. 
Suppose that $(\nabla^i\gamma)(t_0) = 0, 1 \leq i < k$ and that 
$(\nabla^k\gamma)(t_0), (\nabla^{k+1}\gamma)(t_0)$ are linearly independent. 
Then the germ of $\nabla{\mbox{\rm -}}\Tan(\gamma)$ is a frontal with non-degenerate singular point at $(t_0, 0)$ 
and with the singular locus $S(\nabla{\mbox{\rm -}}\Tan(\gamma)) = \{ s = 0\}$. 
\enl

To prove Lemma \ref{non-degenerate-frontal2} we prepare 

\bel
\label{Lemma-}
Let $k \geq 2$. 
Suppose $(\nabla^i\gamma)(t_0) = 0,  1 \leq i < k$ and $(\nabla^k\gamma)(t_0) \not= 0$. 
Then we have: 

{\rm (1)}
For any coordinates of $M$ around $\gamma(t_0)$, 
$\gamma^{(i)}(t_0) = 0, 1 \leq i < k$ and $\gamma^{(k)}(t_0) = (\nabla^k\gamma)(t_0) \not= 0$. 
Moreover we have $\gamma^{(k+1)}(t_0) = (\nabla^{k+1}\gamma)(t_0)$. 

{\rm (2)} 
Set 
$$
u(t) = \frac{1}{k(t - t_0)^{k-1}}\gamma'(t). 
$$
Then $u$ is a $C^\infty$ vector field along $\gamma$ on a neighbourhood of $t_0$. 
The curve $\gamma$ is directed on a neighbourhood of $t_0$ by the frame $u$. 

{\rm (3)}
For any frame $u(t)$ of the directed curve $\gamma$ around $t_0$, and for any $\ell \geq 0$, 
$$
(\nabla^k\gamma)(t_0), \ (\nabla^{k+1}\gamma)(t_0), \ \dots, \ (\nabla^{k+\ell}\gamma)(t_0)
$$ 
are linearly independent if and only if 
$$
u(t_0), \ (\nabla^\gamma_{\pa/\pa t}u)(t_0), \ \dots, \ ((\nabla^\gamma_{\pa/\pa t})^\ell u)(t_0)
$$ 
are linearly independent. 
In particular, for the frame in {\rm (2)}, we have 
$$
u(t_0) = \frac{1}{k!}(\nabla^k\gamma)(t_0), \
(\nabla u)(t_0) 
= \frac{1}{k\cdot k!}(\nabla^{k+1}\gamma)(t_0), \ \dots, \ 
(\nabla^{\ell} u)(t_0) 
= \frac{\ell !}{k\cdot (k+\ell - 1)!}(\nabla^{k+\ell}\gamma)(t_0). 
$$
where $\nabla^i u = (\nabla^\gamma_{\pa/\pa t})^i u$. 
\enl

\Proof 
(1) 
Let $k = 2$. Then $\gamma'(t_0) = (\nabla \gamma)(t_0) = 0$. By Lemma 2.4 of \cite{IY}, 
we have $\gamma''(t_0) = (\nabla^2\gamma)(t_0) \not= 0, \gamma'''(t_0) = (\nabla^3\gamma)(t_0)$. 
Let $k \geq 3$. Then 
$(\nabla^k\gamma)^\lambda$ is a sum of $(\gamma^{(k)})^\lambda$ 
and a polynomial of $\Gamma^\lambda_{\mu\nu}$, their partial derivatives 
and $\gamma^{(i)}, i < k$, each monomial of which 
contains a $\gamma^{(i)}$ with $i \leq k-2$ (cf. Lemma 2.4 of \cite{IY}). 
Thus we have $\gamma^{(i)}(t_0) = (\nabla^i\gamma)(t_0) = 0,  1 \leq i < k$. Moreover we have 
$0 \not= (\nabla^k\gamma)(t_0) = \gamma^{(k)}(t_0)$ and $(\nabla^{k+1}\gamma)(t_0)  = \gamma^{(k+1)}(t_0)$. 
\\
(2) is clear. 
\\
(3) 
We have that $c(t)u(t) = \gamma'(t)$ 
for some function $c(t)$. If $k \geq 2$, then $c(t_0) = 0$. 
By operating $\nabla^\gamma_{\pa/\pa t}$ to both sides of $c(t)u(t) = \gamma'(t)$, 
we have 
$$
c'(t)u(t) + c(t)(\nabla^\gamma_{\pa/\pa t}\,u)(t) = (\nabla^2\gamma)(t). 
$$
If $k \geq 3$, then $c(t_0) = 0, c'(t_0) = 0$. 
In general we have 
$$
c(t_0) = c'(t_0) = \dots = c^{(k-2)}(t_0) = 0, c^{(k-1)}(t_0) \not= 0, 
$$
and 
$$
\begin{array}{ccc}
c^{(k-1)}(t)u(t) + (k-1)c^{(k-2)}(t)(\nabla u)(t) + { }_{k-1}C_2 c^{(k-3)}(t)(\nabla^2 u)(t) + \cdots & = & (\nabla^{k}\gamma)(t)
\\
c^{(k)}(t)u(t) + kc^{(k-1)}(t)(\nabla u)(t) + { }_{k}C_2 c^{(k-2)}(t)(\nabla^2 u)(t) + \cdots & = & (\nabla^{k+1}\gamma)(t)
\\
\vdots & & \vdots
\\
c^{(k+\ell -1)}(t)u(t) + \cdots + { }_{k+\ell - 1}C_{k-1} c^{(k-1)}(t)(\nabla^\ell u)(t) + \cdots & = & 
(\nabla^{k+\ell}\gamma)(t). 
\end{array}
$$
Evaluating at $t_0$, we have the result. 
\QED

\

\noindent
{\it Proof of Lemma \ref{non-degenerate-frontal2}.} 
\ 
The case $k = 1$ is proved in Lemma 3.1 of \cite{IY}. 
Therefor we suppose $k \geq 2$. Let $u(t)$ be a frame around $t_0$ of the directed curve $\gamma$ 
and $c(t)u(t) = \gamma'(t)$, $u(t_0) \not= 0$. 
(For instance $c(t) = k(t - t_0)^{k-1}$). 
Since $f(t, s) = \gamma(t) + s u(t) + \frac{1}{2}s^2 h(\gamma(t), u(t), s)$, we have 
\begin{align*}
\frac{\pa f}{\pa t} & =  \gamma' + su' + 
\frac{1}{2}s^2\, (\gamma')^\mu\frac{\pa h}{\pa x^\mu}(\gamma, u, s) + 
\frac{1}{2}s^2\, (u')^\nu\frac{\pa h}{\pa v^\nu}(\gamma, u, s), 
\\
\frac{\pa f}{\pa s} & =  u + s\, h(\gamma, u, s) + \frac{1}{2}s^2\, \frac{\pa h}{\pa s}(\gamma, u, s). 
\end{align*}
Then we see that $S(f) \supseteq \{ s = 0 \}$ and the kernel field of $f_*$ along $\{ s = 0\}$ is given by 
$\eta = \frac{\pa}{\pa t} - c(t)\frac{\pa}{\pa s}$. 
Let $s \not= 0$. Then
\begin{equation*}
\begin{split}
\frac{1}{s}(\frac{\pa f}{\pa t} - c(t)\frac{\pa f}{\pa s}) = 
& \ 
u' + \frac{1}{2}s\, (\gamma')^\mu\frac{\pa h}{\pa x^\mu}(\gamma, u, s) 
+ \frac{1}{2}s\, (u')^\nu\frac{\pa h}{\pa v^\nu}(\gamma, u, s) 
\\
& 
\hspace{0.5truecm} 
 - c(t) h(\gamma, u, s) - \frac{1}{2}s c(t) \frac{\pa h}{\pa s}(\gamma, u, s). 
\end{split}
\end{equation*}
We define $F(t, s)$ by the right hand side. Then  
$F(t, s) = \frac{1}{s}(\frac{\pa f}{\pa t} - c(t)\frac{\pa f}{\pa s})$ if 
$s \not= 0$. Moreover $F$ is $C^\infty$ also on $s = 0$ and 
$$
F(t, 0) = u'(t) - c(t) h(\gamma(t), u(t), 0). 
$$
By Lemmas 2.1 and 2.2 of \cite{IY}, 
$$
F(t, 0) = u'(t) + c(t)\Gamma^\lambda_{\mu\nu}(\gamma(t))\, (u(t))^\mu(u(t))^\nu 
= (\nabla^\gamma_{\pa/\pa t}u)(t). 
$$
By Lemma \ref{Lemma-} (3), 
if 
$(\nabla^k \gamma)(t_0), (\nabla^{k+1} \gamma)(t_0)$ are linearly independent, then 
$\frac{\pa f}{\pa s}(t, s)$ and $F(t, s)$ are linearly independent around $(t_0, 0)$.  
Moreover they satisfies 
$$
(\frac{\pa f}{\pa t} \wedge \frac{\pa f}{\pa s})(t, s) = -s (\frac{\pa f}{\pa s} \wedge F)(t, s). 
$$
Therefore we see that $\frac{\pa f}{\pa s}(t, s)$ and $F(t, s)$ define an integral lifting of $f$, $f$ is frontal with 
non-degenerate singular point at $(t_0, 0)$, and that 
$S(f) = \{ s = 0\}$. 
\QED

\section{Swallowtails and open swallowtails}
\label{Swallowtail and open swallowtail}

Let $g : (\R^n, p) \to (\R^\ell, q)$ be a map-germ. 
A map germ $f : (\R^n, p) \to \R^{\ell+r}$ is called an {\it opening} of $g$ 
if $f$ is of form $f = (g, h_1, \dots, h_r)$ for some functions $h_1, \dots, h_r : (\R^n, p) \to \R$ 
satisfying 
$$
dh_i = \sum_{j=1}^\ell a_{ij}dg_j, 
$$
for some functions $a_{ij} : (\R^n, p) \to \R, (1 \leq i \leq r, 1 \leq j \leq \ell)$ 
(see for example \cite{Ishikawa4}). 
If $\ell = n$, then the condition on $h$ is equivalent to that 
$f$ is frontal associated with an integral lifting $\widetilde{f} : (\R^n, p) 
\to \Gr(n, T\R^{n+r})$ having Grassmannian coordinates $(a_{ij})$  
such that $\widetilde{f}(p)$ projects isomorphically to 
$T_{g(p)}\R^n$ by the projection $\R^{n+r} = \R^n\times \R^r \to \R^n$. 

Based on results in \cite{KRSUY} and \cite{Ishikawa4}, 
we summarize the characterization results on openings of the Whitney's cusp map-germ: 

\bet
\label{characterization3}
Let $f : (\R^2, p) \to M^m, m \geq 2$ be a germ of frontal with a non-degenerate singular point at $p$, 
$V_1, V_2 : (\R^2, p) \to TM$ an associated frame with $\widetilde{f}$ with $V_2(p) \not\in f_*(T_p\R^2)$, 
and $\eta : (\R^2, p) \to T\R^2$ an extension of a kernel field along of $f_*$. 
Let $c : (\R, t_0) \to (\R^2, p)$ be a parametrization of the singular locus of $f$. 
Set $\gamma = f\circ c : (\R, t_0) \to M$. 
Suppose $(\nabla\gamma)(t_0) = 0$ and $(\nabla^2\gamma)(t_0) \not= 0$. 
Then $f$ 
is diffeomorphic to an opening of Whitney's cusp, 
the germ defined by $(u, t) \mapsto (u, t^3 + ut)$. Moreover we have 

{\rm (0)} Let $m = 2$. Then $f$ is diffeomorphic to Whitney's cusp. 

{\rm (1)} Let $m = 3$. Then $f$ is diffeomorphic to the swallowtail if and only if
$$
V_1(c(t_0)), \ V_2(c(t_0)), \ (\nabla^f_\eta V_2)(c(t_0))
$$ 
are linearly independent in $T_{f(p)}M$. 

{\rm (2)} Let $m \geq 4$. Then $f$ is diffeomorphic to the open swallowtail if and only if
$$
(V_1\circ c)(t_0), \ (V_2\circ c)(t_0), \ ((\nabla^f_\eta V_2)\circ c)(t_0), \ 
(\nabla^\gamma_{\pa/\pa t}((\nabla^f_\eta V_2)\circ c))(t_0)
$$ 
are linearly independent in $T_{f(p)}M$.  

Here $\nabla^f_\eta$ means the covariant derivative by a vector field $\eta$ along a mapping 
$f$ (see \cite{IY}\cite{IY-pre}). 
\ent

\Proof
The assertion (0) follows from Whitney's theorem (also see \cite{Whitney}\cite{SUY}\cite{Saji}). 
(1) follows from Proposition 1.3 of \cite{KRSUY}. 
In general cases $m \geq 2$, we see that 
there exists a submersion $\pi : (M, f(p)) \to (\R^2, 0)$ such that $\pi^{-1}(0)$ is 
transverse to $\widetilde{f}(0) \subset T_{f(p)}M$, $\pi\circ f$ satisfies the same condition with $f$, 
namely, that $\pi\circ f$ is a frontal with the non-degenerate singular point at $p$ and with the same singular locus with 
$f$ and $\eta(c(t_0))$ and $c'(t_0)$ are linearly independent, but $m = 2$. Thus by the assertion (0), 
the map-germ $\pi\circ f$ is diffeomorphic to the Whitney's cusp. 
Moreover we see $f$ is an opening of Whitney's cusp because $f$ is frontal. 

Let $f(u, t) = (u, t^3 + ut, h_1(u, t), \dots, h_r(u, t)), m = 2 + r$ and 
$dh_i = a_idu + b_id(t^3 + ut) = (a_i + tb_i)du + (3t^2 + u)b_i dt, $ 
for some functions $a_i = a_i(u, t), b_i = b_i(u, t), 1 \leq i \leq r$. 
Then we have 
$$
\frac{\pa f}{\pa u} = (1, t, a_1 + t b_1, \dots, a_r + t b_r), 
\quad \frac{\pa f}{\pa t} = (0, 3t^2 + u, (3t^2 + u)b_1, \dots, (3t^2 + u)b_r), 
$$
a frame $V_1 = \frac{\pa f}{\pa u}, V_2 = \frac{1}{3t^2 + u}\frac{\pa f}{\pa t} = (0, 1, b_1, \dots, b_r)$
of the frontal $f$, and a kernel field $\eta = \frac{\pa}{\pa t}$ of $f_*$. 
We have
\begin{align*}
V_1(0, 0) = (1, 0, a_1(0, 0), \dots, a_r(0, 0)), \quad V_2(0, 0) = (0, 1, b_1(0, 0), \dots, b_r(0, 0)), 
\\
(\nabla^f_\eta V_2)(0, 0) = (0, 0, \frac{\pa b_1}{\pa t}(0, 0), \dots, \frac{\pa b_r}{\pa t}(0, 0)). 
\end{align*}
Let $c(t) = (- 3t^2, t)$. Then $\gamma(t) = f(c(t)) = (-3t^2, -2t^3, h_1(-3t^2, t), \dots, h_r(-3t^2, t))$ and 
$$
\nabla^f_\eta V_2(c(t)) = (0, 0, \frac{\pa b_1}{\pa t}(c(t)), \dots, \frac{\pa b_r}{\pa t}(c(t))). 
$$
Then we have 
$$
\nabla^\gamma_{\pa/\pa t}((\nabla^f_\eta V_2)\circ c)\vert_{t=0} = 
(0, 0, \frac{\pa^2 b_1}{\pa t^2}(0, 0), \dots, \frac{\pa^2 b_r}{\pa t^2}(0, 0)). 
$$
Thus the condition of (2) is equivalent, in our case, 
to that $f$ is a versal opening of $\pi\circ f$ and 
then we see $f$ is diffeomorphic to the open swallowtail 
(see Proposition 6.8 (3) $\ell = 3$ of \cite{Ishikawa4}). 
Thus we have the characterization (2). 
\QED

\

\noindent
{\it Proof of Theorem \ref{characterization-theorem}.} 
\ 
Theorem \ref{characterization-theorem} (1) is proved in \cite{IY} in regular case (\S 7 of \cite{IY}). 
Suppose that $\gamma : I \to M$ is not an immersion at $t_0$, $\gamma'(t_0) = 0$, but $\gamma''(t_0) \not= 0$. 
Let $c(t)u(t) = \gamma'(t), u(t_0) \not= 0$. 
Then the $\nabla$-tangent surface is defined by $f(t, s) = \varphi(\gamma(t), u(t), s)$ using the 
geodesics $\varphi(x, v, s)$ on $TM$. 
Then we have the frame 
$$
V_1(t, s) = \frac{\pa f}{\pa s}(t, s), \quad 
V_2(t, s) = F(t,s) = \frac{1}{s}(\frac{\pa f}{\pa t} - c(t)\frac{\pa f}{\pa s}). 
$$
We set $\eta = \frac{\pa}{\pa t} - c(t)\frac{\pa}{\pa s}$. 
Then, by Lemma 5.1 of \cite{IY}, we have 
$(\nabla^f_\eta F)(t, 0) = ({\nabla^\gamma_{\pa/\pa t}}^2 \, u)(t)$. 
Therefore we have 
$\nabla^\gamma_{\pa/\pa t}((\nabla^f_\eta F)(t, 0)) = ({\nabla^\gamma_{\pa/\pa t}}^3 \, u)(t)$. 
Now, by Lemma \ref{Lemma-}, 
$$
V_1(t_0, 0), V_2(t_0, 0), (\nabla^f_\eta F)(t_0, 0)
$$ 
are linearly independent if and only if 
$u(t_0),  (\nabla^\gamma_{\pa/\pa t} u)(t_0),  ({\nabla^\gamma_{\pa/\pa t}}^2 \, u)(t_0)$ 
are linearly independent, and the condition is equivalent to that 
$(\nabla^2\gamma)(t_0), (\nabla^3\gamma)(t_0), (\nabla^4\gamma)(t_0)$ are linearly independent. 
Then in the case $m = 3$, by Theorem \ref{characterization3} (1), 
we have Theorem \ref{characterization-theorem} (1) for non-regular case as well. 

Let $m \geq 4$. Then 
$V_1(t_0, 0), V_2(t_0, 0), (\nabla^f_\eta F)(t_0, 0), \nabla^\gamma_{\pa/\pa t}((\nabla^f_\eta F)(t, 0))\vert_{t = t_0}$ 
are linearly independent if and only if 
$u(t_0),  (\nabla^\gamma_{\pa/\pa t} u)(t_0),  ({\nabla^\gamma_{\pa/\pa t}}^2 \, u)(t_0), 
({\nabla^\gamma_{\pa/\pa t}}^3 \, u)(t_0)$ are linearly independent, and the condition is equivalent to that 
$(\nabla^2\gamma)(t_0), (\nabla^3\gamma)(t_0), (\nabla^4\gamma)(t_0), (\nabla^5\gamma)(t_0)$ 
are linearly independent. 
By Theorem \ref{characterization3} (2), we have Theorem \ref{characterization-theorem} (2).

\section{Perturbations of directed curves} 
\label{Proof of the genericity theorem}

To treat directed curves (see \S \ref{Tangent surface to directed curve}), 
we consider $PTM = \Gr(1, TM)$ with the natural projection 
$\pi : PTM \to M$ and 
the tautological subbundle $D \subset TPTM$ 
on the tangent bundle of $PTM$: 
For any $(x, \ell) \in PTM$ and for any $v \in T_{(x, \ell)}PTM$, 
$v \in D_{(x, \ell)}$ if and only if $\pi_*(v) \in \ell \subset T_xM$. 
A curve $\widetilde{\gamma} : I \to PTM$ is called {\it integral} 
if $\widetilde{\gamma}_*(\pa/\pa t) \in D_{\widetilde{\gamma}(t)}$, 
for any $t \in I$. Recall that $\gamma = \pi\circ \widetilde{\gamma}$ 
with the lifting $\widetilde{\gamma}$ is called a {\it directed curve}. 

Let $u : I \to TM$ be a vector field along a curve $\gamma : I \to M$. 
For $t_0 \in I$, we set 
$$
b_i := 
\inf\left\{ k \left\vert \ 
\rank\left( u(t_0), (\nabla^\gamma_{\pa/\pa t} u)(t_0), \dots, ((\nabla^\gamma_{\pa/\pa t})^{k-1} u)(t_0) \right) 
= i \right.\right\}. 
$$
We have $1 \leq b_1 < b_2 < \cdots < b_m$, if each $b_i < \infty$. 
Then we call the strictly increasing sequence $(b_1, b_2, \dots, b_m)$ of 
natural numbers the {\it $\nabla$-type} of $u$ at $t_0$. 

Moreover the {\it $\nabla$-type} of a curve $\gamma : I \to M$ itself is defined by the 
$\nabla$-type of the velocity vector field $\gamma' : I \to TM$ along $\gamma$. 

\

Let $\gamma : (\R, t_0) \to M$ be a germ of directed curve 
with an integral lifting $\widetilde{\gamma} : (\R, t_0)
\to PTM$ generated by a frame $u : (\R, t_0) \to TM$, $u(t_0) \not= 0$. 
Then $b_1 = 1$ for $u$, since $u(t_0) \not= 0$.

Then we have

\bep
\label{genericity2}
Let $M$ be a manifold of dimension $m$ with an affine connection $\nabla$. 
Then there exists an open dense subset ${\mathcal O}$ in the space of $C^\infty$ integral curves 
$I \to PTM$ with Whitney $C^\infty$ topology 
such that for any $\widetilde{\gamma} \in {\mathcal O}$ and for any $t_0 \in I$, 
$\gamma = \pi\circ \widetilde{\gamma} : I \to M$ is of $\nabla$-type 
$$
(1, 2, 3), (1, 2, 4),  {\mbox{\rm \ or \ }} (2, 3, 4), 
$$
if $m = \dim(M) = 3$, and 
$$
(1, 2, 3, 4, \dots, m-1, m), (1, 2, 3, 4, \dots, m-1, m+1), {\mbox{\rm \ or \ }} (2, 3, 4, 5, \dots, m, m+1). 
$$
if $m \geq 4$, at $t_0$. 
\enp

\

To show Proposition \ref{genericity2}, we use the following generalization 
of Lemma \ref{Lemma-} (3):  

\bel
\label{c, u}
If $\nabla$-type of $u$ is $(1, b_2, \dots, b_m)$ and the order of $c$ at $t_0$ is $\ell$, 
that is, $c(t_0) = \cdots = c^{(\ell - 1)}(t_0) = 0, c^{(\ell)}(t_0) \not= 0$, then 
$\gamma$ is of $\nabla$-type $(a_1, a_2, \dots, a_m) = (1+\ell, b_2+\ell, \dots, b_m+\ell)$. 
\enl

\Proof
By taking covariant derivative $\nabla$ $\ell$-times of the both sides of $c(t)u(t) = \gamma'(t)$, we have 
$(\nabla\gamma)(t_0) = \cdots = (\nabla^\ell\gamma)(t_0) = 0, (\nabla^{\ell+1}\gamma)(t_0) 
= c^{(\ell)}(t_0)u(t_0) \not= 0$. 
Then 
$$
\rank\left( (\nabla\gamma)(t_0), \dots, (\nabla^\ell\gamma)(t_0), (\nabla^{\ell+1}\gamma)(t_0)\right) 
= \rank \left( c^{(\ell)}(t_0)u(t_0)\right) = 1, 
$$ 
and we have $a_1 = 1+\ell$. 
Moreover we have 
\begin{equation*}
\begin{split}
\rank\left( (\nabla\gamma)(t_0), \dots, (\nabla^\ell\gamma)(t_0), (\nabla^{\ell+1}\gamma)(t_0), (\nabla^{\ell+2}\gamma)(t_0)\right) 
= \rank\left( (\nabla^{\ell+1}\gamma)(t_0), (\nabla^{\ell+2}\gamma)(t_0)\right) 
\\ 
= \rank\left( c^{(\ell)}(t_0)u(t_0), c^{(\ell+1)}(t_0)u(t_0) + (\ell+1)c^{(\ell)}(t_0)\nabla u(t_0)\right) 
= \rank\left( u(t_0), \nabla u(t_0)\right). 
\end{split}
\end{equation*}
In general, we have inductively 
$$
\rank\left( (\nabla\gamma)(t_0), (\nabla^2\gamma)(t_0), \dots, (\nabla^k\gamma)(t_0) \right) 
= 
\rank\left( u(t_0), (\nabla u)(t_0), \dots, (\nabla^{k-\ell-1} u)(t_0) \right), 
$$
for any $k \geq 1 + \ell$. 
Therefore we have $a_i = b_i + \ell, 1 \leq i \leq m$. 
\QED

\

We need also the following lemma on local perturbations of integral curves. 

\bel
\label{perturbation}
Let $a < t_1 < t_2 < b$ and $\widetilde{\gamma}, \widetilde{\alpha} : (a, b) \to PT\R^m$ be 
integral curves. Then there exists an integral curve $\widetilde{\beta} : (a, b) \to 
PT\R^m$ such that $\widetilde{\beta}(t) = \widetilde{\alpha}(t), a < t \leq t_1$ and 
$\widetilde{\beta}(t) = \widetilde{\gamma}(t), t_2 \leq t < b$. If $\widetilde{\alpha}$ 
is sufficiently close to $\widetilde{\gamma}$ on $[t_1, t_2]$ in Whitney $C^\infty$ topology, 
then $\widetilde{\beta}$ can be taken to be close to $\widetilde{\gamma}$ 
on $(a, b)$ in Whitney $C^\infty$ topology. 
\enl

\Proof
Let $x = (x^\lambda)$ be a system of coordinates of $\R^m$ and $(x, \xi) = (x^\lambda, \xi_\lambda)$ be 
the associated system of coordinates of $T\R^m$. 
Let $(x\circ \widetilde{\gamma})'(t) = c(t)u(t)$, $(x\circ \widetilde{\alpha})'(t) = e(t)v(t)$, for 
some $c, e : (a, b) \to \R$ and $u, v : (a, b) \to \R^m \setminus \{ 0\}$. 
Then we take a function $f : (a, b) \to \R$ and $w : (a, b) \to \R^m$ such that 
$f(t) = e(t)$ on $(a, t_1]$, $f(t) = c(t)$ on $[t_2, b)$, 
$w(t) = v(t)$ on $(a, t_1]$, $w(t) = u(t)$ on $[t_2, b)$ and 
$$
\dint{t_1}{t_2} f(t) w(t) dt = (x\circ\widetilde{\gamma})(t_1) - (x\circ\widetilde{\alpha})(t_1) + \dint{t_1}{t_2} c(t) u(t) dt. 
$$
Then we have the required $\widetilde{\beta}$ by $(\xi\circ \widetilde{\beta})(t) = w(t)$ and 
$$
(x\circ \widetilde{\beta}) (t) = (x\circ\widetilde{\alpha})(t_1) + \dint{t_1}{t} f(t) w(t) dt, \quad (a < t < b). 
$$
\QED

%

\

\noindent
{\it Proof of Proposition \ref{genericity2}.} 
\ 
Let $\widetilde{\gamma} : (\R, t_0) \to PTM$ be a germ of integral curve with 
$\gamma = \pi\circ\widetilde{\gamma}$. 
Let $c(t)u(t) = \gamma'(t)$ for some frame $u : (\R, t_0) \to TM$ 
along $\gamma$, $u(t_0) \not= 0$, and for some function $c : (\R, t_0) \to \R$. 
Note that $\widetilde\gamma$ is determined by the frame $u$. 
The frame $u$ is determined up to the multiplication of functions $b(t)$ with $b(t_0) \not= 0$. 
Given the initial point $q = \gamma(t_0)$, the pair $(u, c)$ determines the directed curve $\gamma$ uniquely. 
Moreover $(\nabla^k u)(t_0) = u^{(k)}(t_0) + Q$, by a polynomial $Q$ 
of $u^{(i)}(t_0), c^{(i)}(t_0), 0 \leq i < k$ and 
$(\pa^\alpha \Gamma^\lambda_{\mu\nu} / \pa x^\alpha)(q), \vert\alpha\vert \leq k - 1$. 
In particular $(\nabla^k u)(t_0)$ depends only on $k$-jet of $(c, u)$ and just on 
the position $q = \gamma(t_0)$. 

Let us consider the $r$-jet bundle $J^r(I, \R\times(TM \setminus \zeta))$ over 
$I\times \R\times(TM \setminus \zeta)$, where $\zeta$ is the zero-section. 
For the projection $I\times \R\times(TM \setminus \zeta) \to I \times M$, take the fibre 
$J^r(I, \R\times(TM \setminus \zeta))_{(t_0, q)}$ 
over a $(t_0, q) \in I\times M$, and consider the set 
\begin{align*}
S_{\nabla} := \ &\{ j^r(c, u)(t_0)  \mid \ u(t_0), (\nabla u)(t_0), 
\dots, (\nabla^{m-1} u)(t_0) 
{\mbox{\rm \ are linearly dependent }} 
\\
& 
\quad {\mbox{\rm and }} 
u(t_0), (\nabla u)(t_0), \dots, (\nabla^{m-2} u)(t_0), (\nabla^{m} u)(t_0) 
{\mbox{\rm \ are linearly dependent}}\}. 
\\
S'_{\nabla} := \ &\{  j^r(c, u)(t_0) \mid \ 
u(t_0), (\nabla u)(t_0), 
\dots, (\nabla^{m-1} u)(t_0) 
{\mbox{\rm \ are linearly dependent}} 
\\
& 
\quad {\mbox{\rm and }} c(t_0) = 0
\}
\\
S''_{\nabla} := \ &\{  j^r(c, u)(t_0) \mid c(t_0) = c'(t_0) = 0 \}. 
\end{align*}
Then, for any but fixed system of local coordinates around $q$ of $M$, 
$S_{\nabla}, S'_{\nabla}, S''_{\nabla}$ are 
algebraic sets of codimension $\geq 2$. 
Let $S_{\nabla}(I, M), S'_{\nabla}(I, M), S''_{\nabla}(I, M)$ be the corresponding 
subbundle of $J^r(I, \R\times(TM \setminus \zeta))$ over $I\times M$. 
For any subinterval $J \subset I$, 
we set  
\begin{align*}
{\widetilde{\mathcal O}}_J :=  
\{ (c, u) : I \to \R \times (TM \setminus \zeta) \mid & \ 
j^r(c, u) : I \to J^r(I, \R\times(TM\setminus\zeta)) 
\\
& 
{\mbox{\rm \ is transverse to }} 
S_{\nabla}(I, M), 
S'_{\nabla}(I, M), S''_{\nabla}(I, M)
{\mbox{\rm \ over \ }} J
\}. 
\end{align*}
Then $\widetilde{{\mathcal O}} = \widetilde{{\mathcal O}}_I$ 
is open dense in Whitney $C^\infty$ topology. 
Let $(c, u) \in \widetilde{{\mathcal O}}$ and $t_0 \in I$. 
Then $j^r(c, u)(t_0) \not \in S_\nabla \cup S'_{\nabla} \cup S''_{\nabla}$. 
Since $j^r(c, u)(t_0) \not\in S_{\nabla}$, 
we have that the $\nabla$-type of $u$ is $(1, 2, \dots, m-1, m)$ or $(1, 2, \dots, m-1, m+1)$. 
Since $j^r(c, u)(t_0) \not\in S'_{\nabla}$, if $c(t_0) = 0$ then 
$\nabla$-type of $u$ must be $(1, 2, \dots, m-1, m)$. 
On the other hand, since $j^r(c, u)(t_0) \not\in S''_{\nabla}$, 
we have that $c(t_0) \not= 0$ or $c(t_0) = 0, c'(t_0) \not= 0$, i.e. 
the order of $c$ at $t_0$ is $0$ or $1$. 
We set 
$$
{\mathcal O}_J := \{ \widetilde{\gamma} : I \to PTM {\mbox{\rm \ integral\ }} \mid \exists (c, u) \in {\widetilde{\mathcal O}}_J, \ \widetilde{\gamma}(t) = \langle u(t)\rangle_{\R}, (\pi\circ\widetilde{\gamma})'(t) = c(t)u(t) \}. 
$$
We will show, for any compact subinterval $J \subset I$, that ${\mathcal O}_J$ is open dense 
and ${\mathcal O} = {\mathcal O}_I$ is open dense in the space of 
integral curves with Whitney $C^\infty$ topology. 

That ${\mathcal O}_J$ and ${\mathcal O}$ are open is clear, since ${\widetilde{\mathcal O}}$ is open. 

We will show ${\mathcal O}_J$ is dense. 
Let $\widetilde{\gamma} : I \to PTM$ be any integral curve and 
${\mathcal I}$ be any open neighbourhood of $\widetilde{\gamma}$. 
We will show ${\mathcal O}_J \cap {\mathcal I} \not= \emptyset$. 
Set $\gamma = \pi\circ{\widetilde{\gamma}} : I \to M$. 
Take any frame $u$ associated to $\widetilde{\gamma}$. 
Then there exists uniquely $c : I \to \R$ which satisfies $c(t)u(t) = \gamma'(t), t \in I$. 
Take a compact subinterval $J' \subset I$ such that 
$J \subsetneq J'$. 
We approximate $(c, u)$ by some $(e, v) \in {\widetilde{\mathcal O}}_J$ and 
that $(e, v) = (c, u)$ outside of $J'$. 
Then $v$ generates a curve $\rho : I \to PTM, \rho(t) = \langle v(t)\rangle_{\R}$, 
which approximates $\widetilde{\gamma}$, however $\rho$ may not be an integral curve. 
Consider the vector field $(\frac{\pa}{\pa t}, e v)$ along the graph of $\pi\circ\rho$ in $I\times M$. 
Extend $(\frac{\pa}{\pa t}, v(t))$ to a vector field $(\frac{\pa}{\pa t}, V(t, x))$ over $I\times M$ with 
a support contained in $I\times K$ for some compact $K \subset M$. Take $t_0 \in J$. 
Take the integral curve $\alpha : I \to M$ 
of the vector field $(\frac{\pa}{\pa t}, e(t)V(t, x))$ through $(t_0, \alpha(t_0))$. 
Then $\alpha'(t) = e(t) V(t, \alpha(t))$. Define the vector field $w : I \to TM$ over $\alpha$ by 
$w(t) = V(t, \alpha(t))$. Then we have $\alpha'(t) = e(t)w(t)$. 
If we choose $(e, v)$ sufficiently close to $(c, u)$, 
then $w(t) \not= 0$ and $(e, w) \in {\widetilde{\mathcal O}}_J$. 
However the integral curve $\widetilde{\alpha}$ defined by $w$ may not belong to ${\mathcal I}$,  
which is an open set for Whitney $C^\infty$ topology. 
Further we modify the perturbation $(e, v)$ over $J' \setminus J$ 
and the extension $V$ over $(J' \setminus J)\times M$ to obtain 
an integral curve $\widetilde{\beta}$ such that $\widetilde{\beta} = \widetilde{\alpha}$ on $J$ and 
$\widetilde{\beta} = \widetilde{\gamma}$ outside of $J'$, using the method of Lemma \ref{perturbation}. 
Then the integral curve $\widetilde{\beta}$ 
approximates $\widetilde{\gamma}$ and belongs to ${\mathcal I}$, while  
$(e, w) \in \widetilde{\mathcal O}_J$. 
Since 
$\widetilde{\beta}(t) = \langle w(t)\rangle_{\R}$ and $(\pi\circ\widetilde{\beta})'(t) = e(t)w(t)$, 
we have $\widetilde{\beta} \in {\mathcal O}_J \cap {\mathcal I}$. 
Thus we have seen that ${\mathcal O}_J$ is dense, for any compact subinterval $J \subset I$. 

Since ${\mathcal O} = \cap_{J \subset I} {\mathcal O}_J$, the intersection over 
compact subintervals $J \subset I$, we have that 
${\mathcal O}$ is residual, and therefore that ${\mathcal O}$ is dense in Whitney $C^\infty$ topology 
\cite{GG}. 

Thus we have that ${\mathcal O}$ is open dense in Whitney $C^\infty$ topology. 
Then, using Lemma \ref{c, u}, 
we have the required result. 
\QED

\ber
{\rm 
By the same method as above, we have that 
the codimension of jets of integral curves such that the projections are of $\nabla$-type 
$(a_1, a_2, \dots, a_m)$ is given by 
$$
\ell + \sum_{i=1}^m (b_i - i) =  a_1 - 1 + \sum_{i=2}^m (a_i - a_1 - i +1), 
$$
for any affine connection $\nabla$. Note that 
the codimension is calculated in Theorem 5.6 of \cite{Ishikawa4} in 
the flat case 
(cf. Theorem 5.8, Theorem 3.3 of \cite{Ishikawa4}). 
}
\enr

\

\noindent
{\it Proof of Theorem \ref{genericity-theorem}.}
\ 
We observe that the equation on geodesics 
$$
\dfrac{\pa^2\varphi}{\pa s^2}^\lambda(x, v, s) \ + \ \Gamma^\lambda_{\mu \nu}(\varphi(x, v, s))
\dfrac{\pa \varphi}{\pa s}^\mu(x, v, s)\dfrac{\pa \varphi}{\pa s}^\nu(x, v, s) = 0, 
$$
is symmetric on the indices $\mu, \nu$. Therefore 
the geodesics $\varphi(x, v, s)$ and the tangent surfaces $\nabla{\mbox{\rm -}}\Tan(\gamma)$ 
remain same if the connection $\Gamma^\lambda_{\mu \nu}$ is replaced by the torsion free 
connection $\frac{1}{2}(\Gamma^\lambda_{\mu \nu} + \Gamma^\lambda_{\nu \mu})$, 
in other word, if $\nabla$ is replaced by the torsion free connection $\widetilde{\nabla}$, 
defined by $\widetilde{\nabla}_XY = \nabla_XY - \frac{1}{2}T(X, Y)$. 
Thus we may suppose $\nabla$ is torsion free. Then Theorem \ref{characterization-theorem} 
and Proposition \ref{genericity2} imply Theorem \ref{genericity-theorem}. 
\QED

\

\begin{flushleft}
Goo ISHIKAWA, \\
Department of Mathematics, Hokkaido University, 
Sapporo 060-0810, Japan. \\
e-mail : ishikawa@math.sci.hokudai.ac.jp \\

\

Tatsuya YAMASHITA, \\
Department of Mathematics, Hokkaido University, 
Sapporo 060-0810, Japan. \\
e-mail : tatsuya-y@math.sci.hokudai.ac.jp \\ 

\end{flushleft}

\end{document}